\documentclass[11pt]{article}
\usepackage{amssymb,amsthm,amsmath}
\def\qed{\hfill $\Box$}

\newcommand\pf{\smallbreak\noindent \texttt{Proof}. }

\begin{document}

\newtheorem{thm}{Theorem}[section]
\newtheorem{prop}[thm]{Proposition}
\newtheorem{lem}[thm]{Lemma}
\newtheorem{cor}[thm]{Corollary}
\newtheorem{ex}[thm]{Example}
\renewcommand{\thefootnote}{*}

\title{\bf Relationships between the generalized central series of Leibniz algebras}

\author{\textbf{Aleksandr~A.~Pypka}\\
Oles Honchar Dnipro National University, Dnipro, Ukraine\\
{\small e-mail: sasha.pypka@gmail.com}}
\date{}

\maketitle

\begin{abstract}
The purpose of this article is to show a close relationship between the generalized central series of Leibniz algebras. Some analogues of the classical group-theoretical theorems of Schur and Baer for Leibniz algebras are proved.
\end{abstract}



\thispagestyle{empty}

\section{Introduction}
We start our motivation with some classical group-theoretical results. In 1951, B.H.~Neumann proved~\cite{NB1951} the so-called~\cite{KS2016} Schur's theorem: if $G$ is a group and $G/\zeta(G)$ is finite, then $[G,G]$ is also finite and $|[G,G]|\leqslant t^{t^{2}+1}$, where $t=|G/\zeta(G)|$. Later J.~Wiegold~\cite{WJ1965} and B.A.F.~Wehrfritz~\cite{WB2018} significantly improved this bound. There are many analogues of Schur's theorem for groups (see, for example, \cite{BCKO2013, DKP2015, DKS2018, FGK1995, KOP2016A, KSP2014, KS2013, MA2007, P2014A, P2017}, etc.). Among them there are results related to automorphism groups. Let $G$ be a group and $A$ be a subgroup of the automorphism group $\mathrm{Aut}(G)$ of $G$. Put
\begin{gather*}
C_{G}(A)=\{g\in G|\ \alpha(g)=g\ \mbox{for each }\alpha\in A\},\\
[G,A]=\langle[g,\alpha]=g^{-1}\alpha(g)|\ g\in G,\alpha\in A\rangle.
\end{gather*}
The subgroups $C_{G}(A)$ and $[G,A]$ are called the \textit{$A$-center} and \textit{$A$-commutator subgroup} of $G$. P.~Hegarty proved~\cite{HP1994} that if $A=\mathrm{Aut}(G)$ and $G/C_{G}(A)$ is finite then $[G,A]$ is also finite. M.R.~Dixon, L.A.~Kurdachenko and A.A.~Pypka in \cite{DKP2014} considered the following more general situation: $\mathrm{Inn}(G)\leqslant A$ and $A/\mathrm{Inn}(G)$ is finite, where $\mathrm{Inn}(G)$ is the group of inner automorphisms of $G$. For this case the automorphic analogue of Schur's theorem was obtained: if $G$ is a group and $G/C_{G}(A)$ is finite then $[G,A]$ is also finite. Moreover, the bound for the order of $[G,A]$ in the terms of the orders $|G/C_{G}(A)|$ and $|A/\mathrm{Inn}(G)|$ was obtained. In particular, when $A=\mathrm{Inn}(G)$ or $A=\mathrm{Aut}(G)$ then we just have Schur's and Hegarty's theorems. Other automorphic analogues of Schur's theorem can be found in \cite{CDR2016,KOP2015,KP2014,P2014}.

It is well known that there is a very close relationship between groups and Lie algebras. For Lie algebras, many group-theoretical results are proved, and vice versa. This topic is no exception. In particular, for Lie algebras an analogue of Schur's theorem is well-known (see, for example, \cite{VL1972}). E.L.~Stitzinger and R.M.~Turner proved \cite{ST1999} a Lie analogue of Hegarty's theorem. Let $L$ be a Lie algebra over a field $F$ and $\mathrm{Der}(L)$ be an algebra of derivations of $L$. According to \cite{ST1999} we put
$$H=\bigcap\limits_{\alpha\in\mathrm{Der}(L)}\mathrm{Ker}(\alpha)\ \mbox{ and }\ L^{\ast}=\sum\limits_{\alpha\in\mathrm{Der}(L)}\mathrm{Im}(\alpha).$$
Thus the subalgebras $H$ and $L^{\ast}$ are Lie analogues of $C_{G}(\mathrm{Aut}(G))$ and $[G,\mathrm{Aut}(G)]$. From the main result of \cite{ST1999} follows that if $L/H$ has finite dimension, then $L^{\ast}$ is also finite-dimensional.

Consider now the (left) Leibniz algebras. Let $L$ be an algebra over a field $F$ with the binary operations $+$ and $[\ ,\ ]$. Then $L$ is called a \textit{left Leibniz algebra} if $[[x,y],z]=[x,[y,z]]-[y,[x,z]]$ for all $x,y,z\in L$ \cite{Blo1965,Lo1993}. Note that every Lie algebra $L$ is a Leibniz algebra. Moreover, Lie algebras can be characterized as Leibniz algebras in which $[x,x]=0$ for every $x\in L$.

One of the directions in the development of Leibniz algebra theory is the search for analogues of the results of Lie algebra theory. At the same time, there is a rather significant difference between these types of algebras (see, for example, \cite{CPSY2019,KKPS2017}). Recall some definitions. The \textit{left} (respectively \textit{right}) \textit{center} $\zeta^{l}(L)$ (respectively $\zeta^{r}(L)$) of a Leibniz algebra $L$ is defined by the rule: $\zeta^{l}(L)=\{x\in L|\ [x,y]=0\ \mbox{for each }y\in L\}$ (respectively $\zeta^{r}(L)=\{x\in L|\ [y,x]=0\ \mbox{for each }y\in L\}$).
The left center is an ideal of $L$, but it is not true for the right center. The right center is a subalgebra of $L$, and in general, the left and right centers are different. Moreover, they even may have different dimensions (see example~2.1 in \cite{KOP2016}). The \textit{center} $\zeta(L)$ of $L$ is the intersection of the left and right centers, that is
$$\zeta(L)=\{x\in L|\ [x,y]=0=[y,x]\ \mbox{for each }y\in L\}.$$
Clearly, $\zeta(L)$ is an ideal of $L$. Thus, we can consider the factor-algebra $L/\zeta(L$).

L.A.~Kurdachenko, J.~Otal and A.A.~Pypka proved \cite{KOP2016} the following modification of Schur's theorem: if $L$ is a Leibniz algebra over a field $F$, $\mathrm{codim}_{F}(\zeta^{l}(L))=l$ and $\mathrm{codim}_{F}(\zeta^{r}(L))=r$ are finite, then $\mathrm{dim}_{F}([L,L])\leqslant l(l+r)$. In this regard, the following question arises. Suppose that only $\mathrm{codim}_{F}(\zeta^{l}(L))$ is finite. Is $\mathrm{dim}_{F}([L,L])$ finite? An example~3.1 from \cite{KOP2016} gives a negative answer on this question. However, we have the direct analogue of Schur's theorem for Leibniz algebras: if $L$ is a Leibniz algebra over a field $F$ and $\mathrm{codim}_{F}(\zeta(L))=d$ is finite, then $\mathrm{dim}_{F}([L,L])\leqslant d^{2}$ \cite{KOP2016}.

Taking into account the above arguments, it is natural to consider an analogue of the result of \cite{DKP2014} for Leibniz algebras. First we define the analogues of the $A$-center and $A$-commutator subgroup for Leibniz algebras. Let $L$ be a Leibniz algebra over a field $F$ and $D$ be a subalgebra of $\mathrm{Der}(L)$. Put
$$\mathrm{Ann}_{L}(D)=\bigcap\limits_{\alpha\in D}\mathrm{Ker}(\alpha)\ \mbox{ and }\ [L,D]=\sum\limits_{\alpha\in D}\mathrm{Im}(\alpha).$$

Let $a\in L$. Consider the mapping $\mathrm{l}_{a}:L\rightarrow L$ defined by $\mathrm{l}_{a}(x)=[a,x]$. It is well known that $\mathrm{l}_{a}$ is a derivation of $L$ and the set $\mathrm{Ad}^{l}(L)=\{\mathrm{l}_{a}|\ a\in L\}$ is an ideal of $\mathrm{Der}(L)$ (see, for example, \cite{KKPS2017}).

Suppose that $\mathrm{Ad}^{l}(L)\leqslant D$, then $\mathrm{Ann}_{L}(D)\leqslant\mathrm{Ann}_{L}(\mathrm{Ad}^{l}(L))=\zeta^{r}(L)$. Thus
$$\mathrm{Ann}_{L}(D)\cap\zeta^{l}(L)\leqslant\mathrm{Ann}_{L}(\mathrm{Ad}^{l}(L))\cap\zeta^{l}(L)=\zeta^{r}(L)\cap\zeta^{l}(L)=\zeta(L).$$
In particular, $\mathrm{Ann}_{L}(D)\cap\zeta^{l}(L)$ is an ideal of $L$.

We will say that $\mathrm{A}_{L}(D)=\mathrm{Ann}_{L}(D)\cap\zeta^{l}(L)$ is the \textit{$D$-center} of $L$. Note that $\mathrm{A}_{L}(D)\leqslant\zeta(L)$ and if $D=\mathrm{Ad}^{l}(L)$ then the $D$-center is a usual center of $L$. We will say that $[L,D]$ is the \textit{$D$-derived subalgebra} of $L$. If $D=\mathrm{Ad}^{l}(L)$ then $[L,D]=[L,\mathrm{Ad}^{l}(L)]$. Let $x\in [L,\mathrm{Ad}^{l}(L)]$, then $x=[y,\mathrm{l}_{a}]=\mathrm{l}_{a}(y)=[a,y]$ for all $a,y\in L$. This means that in this case the $D$-derived subalgebra $[L,D]$ is a usual derived subalgebra $[L,L]$ of $L$.

The first main result of this paper is the following

\textbf{Theorem~A}. \textit{Let $L$ be a Leibniz algebra over a field $F$, $D$ be a subalgebra of $\mathrm{Der}(L)$ such that $\mathrm{Ad}^{l}(L)\leqslant D$. Suppose that $\mathrm{dim}_{F}(D/\mathrm{Ad}^{l}(L))=k$ is finite. If $\mathrm{dim}_{F}(L/\mathrm{A}_{L}(D))=t$ is finite, then $\mathrm{dim}_{F}([L,D])\leqslant t(k+t)$.}

In \cite{BR1952}, R.~Baer generalized Schur's theorem in the following way. He proved that if the factor-group $G/\zeta_{k}(G)$ is finite, then $\gamma_{k+1}(G)$ is also finite. For a group $G$ we denote by $\zeta_{k}(G)$ the $k$th term of the upper central series of $G$ and denote by $\gamma_{k}(G)$ the $k$th term of the lower central series of $G$. An automorphic analogue of this result for groups was obtained in \cite{DKP2014}. If $L$ is a Lie algebra such that $L/\zeta_{k}(L)$ is finite-dimensional, I.N.~Stewart~\cite{SI1974} showed that $\gamma_{k+1}(L)$ is also finite-dimensional. Moreover, the upper bound for the dimension $\mathrm{dim}_{F}(\gamma_{k+1}(L))$ was obtained in \cite{KOP2016}. In the same paper, an analogue of Baer's theorem for Leibniz algebras was proved.

Starting with $\mathrm{A}_{L}(D)$ and $[L,D]$, where $\mathrm{Ad}^{l}(L)\leqslant D$, we can define the upper and lower $D$-central series of $L$. First we let $\zeta_{1}(L,D)=\mathrm{A}_{L}(D)$. This allows us to define an ascending series of ideals of $L$, with terms $\zeta_{\nu}(L,D)$, where $\zeta_{\nu+1}(L,D)/\zeta_{\nu}(L,D)=\zeta_{1}(L/\zeta_{\nu}(L,D),D)$. As usual, if $\lambda$ is a limit ordinal, then we let $\zeta_{\lambda}(L,D)=\bigcup_{\mu<\lambda}\zeta_{\mu}(L,D)$. The last term $\zeta_{\infty}(L,D)=\zeta_{\delta}(L,D)$ of this series is called the \textit{upper $D$-hypercenter} of $L$ and the ordinal $\delta$ is called the \textit{$D$-upper central length} of $L$, which we denote by $\mathrm{zl}(L,D)$.

The \textit{lower $D$-central series} of $L$ is the descending series
$$L=\gamma_{1}(L,D)\geqslant\gamma_{2}(L,D)\geqslant\ldots\gamma_{\nu}(L,D)\geqslant\gamma_{\nu+1}(L,D)\geqslant\ldots$$
defined by $\gamma_{2}(L,D)=[L,D]$, and for each ordinal $\nu$ we have that $\gamma_{\nu+1}(L,D)=[\gamma_{\nu}(L,D),D]$. As usual, for the limit ordinal $\lambda$, we define $\gamma_{\lambda}(L,D)=\bigcap_{\mu<\lambda}\gamma_{\mu}(L,D)$.

The second main result of this paper is the following

\textbf{Theorem~B}. \textit{Let $L$ be a Leibniz algebra over a field $F$, $D$ be a subalgebra of $\mathrm{Der}(L)$ such that $\mathrm{Ad}^{l}(L)\leqslant D$ and $\mathrm{dim}_{F}(D/\mathrm{Ad}^{l}(L))=k$ is finite. Let $Z$ be the upper $D$-hypercenter of $L$. Suppose that $\mathrm{zl}(L,D)=m$ and $\mathrm{dim}_{F}(L/Z)=t$ are finite. Then $\gamma_{m+1}(L,D)$ is finite-dimensional, and there exists a function $f$ such that $\mathrm{dim}_{F}(\gamma_{m+1}(L,D))\leqslant f(k,m,t)$.}

\section{Proof of Theorem A}
\pf Since $\mathrm{Ad}^{l}(L)\leqslant D$, $\mathrm{A}_{L}(D)\leqslant\zeta(L)$, which implies that
$$\mathrm{dim}_{F}(L/\zeta(L))\leqslant\mathrm{dim}_{F}(L/\mathrm{A}_{L}(D))=t.$$
Then $K=[L,L]$ is finite-dimensional and $\mathrm{dim}_{F}([L,L])\leqslant t^{2}$~\cite{KOP2016}.

Put $L_{ab}=L/K$. For each $\alpha\in D$ we define the mapping $\alpha_{ab}:L_{ab}\rightarrow L_{ab}$ by the following rule: $\alpha_{ab}(x+K)=\alpha(x)+K$ for every $x\in L$. Let $x,y\in L$, $\lambda\in F$. Since
\begin{gather*}
\alpha_{ab}((x+K)+(y+K))=\\
=\alpha_{ab}((x+y)+K)=\alpha(x+y)+K=\\
=\alpha(x)+\alpha(y)+K=(\alpha(x)+K)+(\alpha(y)+K)=\\
=\alpha_{ab}(x+K)+\alpha_{ab}(y+K)
\end{gather*}
and
\begin{gather*}
\alpha_{ab}(\lambda(x+K))=\alpha_{ab}(\lambda x+K)=\\
=\alpha(\lambda x)+K=\lambda\alpha(x)+K=\\
=\lambda(\alpha(x)+K)=\lambda\alpha_{ab}(x+K),
\end{gather*}
$\alpha_{ab}$ is an endomorphism of $L_{ab}$.

Let again $x,y\in L$. Then
\begin{gather*}
\alpha_{ab}([x+K,y+K])=\alpha_{ab}([x,y]+K)=\\
=\alpha([x,y])+K=[\alpha(x),y]+[x,\alpha(y)]+K=\\
=([\alpha(x),y]+K)+([x,\alpha(y)]+K)=\\
=[\alpha(x)+K,y+K]+[x+K,\alpha(y)+K]=\\
=[\alpha_{ab}(x+K),y+K]+[x+K,\alpha_{ab}(y+K)],
\end{gather*}
which shows that $\alpha_{ab}\in\mathrm{Der}(L_{ab})$.

Let $\alpha\in D$ and consider a mapping $d(\alpha):L_{ab}\rightarrow L_{ab}$ defined by the following rule: $d(\alpha)(x)=[x,\alpha_{ab}]=\alpha_{ab}(x)$, $x\in L_{ab}$. We have
\begin{gather*}
d(\alpha)(x+y)=[x+y,\alpha_{ab}]=\alpha_{ab}(x+y)=\\
=\alpha_{ab}(x)+\alpha_{ab}(y)=[x,\alpha_{ab}]+[y,\alpha_{ab}]=\\
=d(\alpha)(x)+d(\alpha)(y),\\
d(\alpha)(\lambda x)=[\lambda x,\alpha_{ab}]=\alpha_{ab}(\lambda x)=\\
=\lambda\alpha_{ab}(x)=\lambda[x,\alpha_{ab}]=\lambda d(\alpha)(x).
\end{gather*}
In other words, $d(\alpha)$ is an endomorphism of $L_{ab}$.

Furthermore, $\mathrm{Im}(d(\alpha))=[L_{ab},\alpha_{ab}]$, $\mathrm{Ker}(d(\alpha))\geqslant\mathrm{A}_{L_{ab}}(\alpha_{ab})$, which implies that $\mathrm{dim}_{F}(L_{ab}/\mathrm{Ker}(d(\alpha)))\leqslant\mathrm{dim}_{F}(L_{ab}/\mathrm{A}_{L_{ab}}(\alpha_{ab}))$. Thus, we have
$$[L_{ab},\alpha_{ab}]=\mathrm{Im}(d(\alpha))\cong L_{ab}/\mathrm{Ker}(d(\alpha)).$$

If $x\in\mathrm{A}_{L}(\alpha)$, then
$$\alpha_{ab}(x+K)=\alpha(x)+K=K,$$
which shows that $x+K\in\mathrm{A}_{L_{ab}}(\alpha_{ab})$. Thus $(\mathrm{A}_{L}(\alpha)+K)/K\leqslant\mathrm{A}_{L_{ab}}(\alpha_{ab})$. An inclusion $\mathrm{A}_{L}(D)\leqslant\mathrm{A}_{L}(\alpha)$ implies that
$$\mathrm{dim}_{F}(L/\mathrm{A}_{L}(\alpha))\leqslant\mathrm{dim}_{F}(L/\mathrm{A}_{L}(D))=t.$$
Thus $\mathrm{dim}_{F}(L_{ab}/\mathrm{A}_{L_{ab}}(\alpha_{ab}))\leqslant t$ and we obtain that $\mathrm{dim}_{F}([L_{ab},\alpha_{ab}])\leqslant t$ for every $\alpha\in D$.

Let $B=\{\alpha_{1},\ldots,\alpha_{k}\}$ be a basis of $D/\mathrm{Ad}^{l}(L)$ and $\beta\in\mathrm{Ad}^{l}(L)+\alpha$. Then $\beta=\mathrm{l}_{z}+\alpha$ for some $z\in L$. For arbitrary element $y\in L$ we have
$$[y,\beta]=\beta(y)=(\mathrm{l}_{z}+\alpha)(y)=\mathrm{l}_{z}(y)+\alpha(y)=[z,y]+[y,\alpha].$$
It follows that
\begin{gather*}
[y,\beta]+K=[z,y]+[y,\alpha]+K=\\
=([z,y]+K)+([y,\alpha]+K)=\\
=[z+K,y+K]+([y,\alpha]+K).
\end{gather*}
Since $L/K$ is abelian, $[z+K,y+K]=0$, that is $[y,\beta]+K=[y,\alpha]+K$. On the other hand,
$$[x+K,\alpha_{ab}]=\alpha_{ab}(x+K)=\alpha(x)+K=[x,\alpha]+K,$$
which shows that $[L_{ab},\alpha_{ab}]=([L,\alpha]+K)/K$.

Clearly, $[L,D]=\langle[L,\beta]|\ \beta\in D\rangle = \langle[L,\beta]|\ \beta\in\mathrm{Ad}^{l}(L)+\alpha_{j}, 1\leqslant j\leqslant k\rangle$, so that
\begin{gather*}
([L,D]+K)/K=\sum_{\substack{\beta\in\mathrm{Ad}^{l}(L)+\alpha_{j}\\ 1\leqslant j\leqslant k}}([L,\beta]+K)/K=\\
=\sum\limits_{1\leqslant j\leqslant k}([L,\alpha_{j}]+K)/K=\sum\limits_{1\leqslant j\leqslant k}[L_{ab},(\alpha_{j})_{ab}].
\end{gather*}
It follows that $\mathrm{dim}_{F}(([L,D]+K)/K)\leqslant tk$, which means that
$$\mathrm{dim}_{F}([L,D])\leqslant tk+t^{2}=t(k+t).$$\qed

If $D=\mathrm{Ad}^{l}(L)$, then $\mathrm{A}_{L}(D)=\zeta(L)$, $[L,D]=[L,L]$ and we obtain the following direct analogue of Schur's theorem for Leibniz algebras.

\begin{cor}[\cite{KOP2016}]
Let $L$ be a Leibniz algebra over a field $F$. If $L/\zeta(L)$ has finite dimension $t$, then $\mathrm{dim}_{F}([L,L])\leqslant t^{2}$.
\end{cor}

In particular, if $L$ is a Lie algebra, we obtain the following direct analogue of Schur's theorem for Lie algebras.

\begin{cor}[\cite{VL1972}]
Let $L$ be a Lie algebra over a field $F$. If $L/\zeta(L)$ has finite dimension $t$, then $\mathrm{dim}_{F}([L,L])\leqslant t(t+1)/2$.
\end{cor}

Repeating almost word-to-word the proof of Theorem~1 from \cite{ST1999}, we can obtain that if $\mathrm{dim}_{F}(L/\mathrm{A}_{L}(\mathrm{Der}(L)))$ is finite, then $\mathrm{Der}(L)$ is finite-dimensional. Thus, if $D=\mathrm{Der}(L)$, then we have the following direct analogue of Hegarty's theorem for Leibniz algebras.

\begin{cor}
Let $L$ be a Leibniz algebra over a field $F$. If $L/\mathrm{A}_{L}(\mathrm{Der}(L))$ has finite dimension $t$, then $\mathrm{dim}_{F}([L,\mathrm{Der}(L)])\leqslant t(t+1)$.
\end{cor}

In particular, if $L$ is a Lie algebra, we obtain the following direct analogue of Hegarty's theorem for Lie algebras.

\begin{cor}[\cite{ST1999}]
Let $L$ be a Lie algebra over a field $F$. If $L/\mathrm{A}_{L}(\mathrm{Der}(L))$ has finite dimension $t$, then $\mathrm{dim}_{F}([L,\mathrm{Der}(L)])\leqslant t(t+1)/2$.
\end{cor}

\section{Proof of Theorem B}
\pf Let
$$\langle0\rangle=Z_{0}\leqslant Z_{1}\leqslant\ldots\leqslant Z_{m-1}\leqslant Z_{m}=Z$$
be the upper $D$-central series of $L$. We use induction on $m$. If $m=1$, then $Z_{1}=\mathrm{A}_{L}(D)$, $\mathrm{dim}_{F}(L/Z_{1})=t$, and Theorem~A implies that $[L,D]=\gamma_{2}(L,D)$ is finite-dimensional of dimension at most $t(k+t)$.

Suppose inductively that the result is true for some $m-1$ and let $L$ be a Leibniz algebra satisfying the hypotheses of the theorem with $\mathrm{zl}(L,D)=m$. Let $U=L/Z_{1}$. For each $\alpha\in D$ we define the mapping $\alpha_{U}:U\rightarrow U$ by the following rule:
$$\alpha_{U}(x+Z_{1})=\alpha(x)+Z_{1}$$
for every $x\in L$. Since
$$\alpha_{U}((x+Z_{1})+(y+Z_{1}))=\alpha_{U}(x+Z_{1})+\alpha_{U}(y+Z_{1})$$
and
$$\alpha_{U}(\lambda(x+Z_{1}))=\lambda\alpha_{U}(x+Z_{1}),$$
$\alpha_{U}$ is an endomorphism of $U$. Moreover,
$$\alpha_{U}([x+Z_{1},y+Z_{1}])=[\alpha_{U}(x+Z_{1}),y+Z_{1}]+[x+Z_{1},\alpha_{U}(y+Z_{1})].$$
Hence $\alpha_{U}\in\mathrm{Der}(U)$.

Consider the mapping $\eta:D\rightarrow\mathrm{Der}(U)$, defined by the following rule: $\eta(\alpha)=\alpha_{U}$. Obviously $\eta$ is a homomorphism. Let $\alpha\in\mathrm{Ad}^{l}(L)$, that is $\alpha=\mathrm{l}_{a}$ for some $a\in L$. Then $\eta(\alpha)=\alpha_{U}=(\mathrm{l}_{a})_{U}$ and
$$(\mathrm{l}_{a})_{U}(x+Z_{1})=\mathrm{l}_{a}(x)+Z_{1}=[a,x]+Z_{1}=[a+Z_{1},x+Z_{1}]=\mathrm{l}_{a+Z_{1}}(x+Z_{1}),$$
which implies that $\eta(\mathrm{Ad}^{l}(L))=\mathrm{Ad}^{l}(U)$. It follows that $\mathrm{Ad}^{l}(U)=\eta(\mathrm{Ad}^{l}(L))\leqslant\eta(D)$ and $\eta(D)/\mathrm{Ad}^{l}(U)$ is finite-dimensional of dimension at most $k$. Furthermore, the series
$$\langle0\rangle=Z_{1}/Z_{1}\leqslant Z_{2}/Z_{1}\leqslant\ldots\leqslant Z_{m-1}/Z_{1}\leqslant Z_{m}/Z_{1}$$
is the upper $D$-central series of $L/Z_{1}$. Since $\mathrm{zl}(L/Z_{1},\eta(D))=m-1$, the induction hypothesis implies that $\gamma_{m}(L/Z_{1},\eta(D))$ is finite-dimensional and there is a function $\beta(k,m-1,t)$ such that $\mathrm{dim}_{F}(\gamma_{m}(L/Z_{1},\eta(D)))\leqslant\beta(k,m-1,t)$.

We have that $\gamma_{m}(L/Z_{1},D)=(\gamma_{m}(L,D)+Z_{1})/Z_{1}$. Let
$$K/Z_{1}=\gamma_{m}(L/Z_{1},\eta(D))=\gamma_{m}(L/Z_{1},D).$$
We note that $K/Z_{1}$ is finite-dimensional and $\mathrm{dim}_{F}(K/Z_{1})\leqslant\beta(k,m-1,t)=r$. Thus we can apply Theorem~A to $K$. We have that $[K,D]$ is finite-dimensional and $\mathrm{dim}_{F}([K,D])\leqslant r(k+r)$. Finally, since
$$\gamma_{m+1}(L,D)=[\gamma_{m}(L,D),D]\leqslant[K,D],$$
$\gamma_{m+1}(L,D)$ is finite-dimensional of dimension at most $r(k+r)=\beta(k,m,t)$.\qed

The function $\beta(k,m,t)$ is defined recursively by $\beta(k,1,t)=t(k+t)$, and
$$\beta(k,m+1,t)=\beta(k,m,t)(k+\beta(k,m,t)).$$

If $D=\mathrm{Ad}^{l}(L)$, then $\zeta_{k}(L,D)=\zeta_{k}(L)$ and $\gamma_{k}(L,D)=\gamma_{k}(L)$ and we obtain the following direct analogue of Baer's theorem for Leibniz algebras.

\begin{cor}[\cite{KOP2016}]
Let $L$ be a Leibniz algebra over a field $F$. If $L/\zeta_{k}(L)$ has finite dimension $t$, then $\mathrm{dim}_{F}(\gamma_{k+1}(L))\leqslant 2^{k-1}t^{k+1}$.
\end{cor}

In particular, if $L$ is a Lie algebra, we obtain the following direct analogue of Baer's theorem for Lie algebras.

\begin{cor}[\cite{KOP2016},\cite{SI1974}]
Let $L$ be a Lie algebra over a field $F$. If $L/\zeta_{k}(L)$ has finite dimension $t$, then $\mathrm{dim}_{F}(\gamma_{k+1}(L))\leqslant t^{k}(t+1)/2$.
\end{cor}

\end{document}